  \date{\today}   

 \documentclass[11pt,a4paper]{article}
 \usepackage{amsmath,amsfonts,amssymb,amsthm}
 \usepackage{enumerate,color}
\usepackage[applemac]{inputenc}
 \usepackage[polutonikogreek,english]{babel}

\theoremstyle{plain}
\newtheorem{theorem}{Theorem}[section]
\newtheorem{proposition}[theorem]{Proposition}
\newtheorem{corollary}[theorem]{Corollary}
\newtheorem{lemma}[theorem]{Lemma}
\newtheorem{definition}{Definition}[section]

\numberwithin{equation}{section}

\DeclareMathOperator{\Card}{Card}%

\renewcommand{\r}{\mathbb{R}}%
\newcommand{\tq}{\, \big| \, }%

\author{Vladimir Gol'dshtein and Marc Troyanov} 
\title{On the naturality of the exterior differential}
 
\begin{document}
 
\maketitle

\begin{abstract}
We  give sufficient conditions for the naturallity of the exterior differential under Sobolev mappings. In other words we study the validity of the equation  $d\, f^* \alpha = f^*\; d\alpha$ for a smooth form $\alpha$ and a Sobolev map $f$. The main results of the paper are Theorems \ref{th.1} and \ref{th.2}.
\medskip

\end{abstract}

\section{Introduction}

One of the main properties of calculus with differential forms is the \emph{naturallity} of the exterior derivative, that is the fact that
for any smooth map $f : U \to \r^n$, where $U$ is a bounded domain in $\r^m$, and any smooth differential form $\alpha$ in $\r^n$,
we have 
\begin{equation}\label{eq.nat}
 df^*\alpha  = f^*d\alpha.
\end{equation}
Note that this equation is just an avatar of the chain rule; its proof can be found in any textbook on differential forms. 

\medskip

For applications in the calculus of variation, non linear elasticity or geometric analysis, it is important to extend this result to non smooth situations. If the map $f$ is smooth and $\alpha$ is   a Sobolev differential form, then 
the pull back $f^*\alpha$ is also a locally  Sobolev differential form and the naturality (\ref{eq.nat}) can be proved by  standard arguments. If both the differential form $\alpha$ and the map $f$ belong to  $W^{1,1}_{loc}$, then the problem is
not well posed and it is not clear under what conditions, should the equation (\ref{eq.nat}) make sense and be proved.

\bigskip

If the differential form $\alpha$ is smooth, then the situation is better and it is our goal in this paper to give sufficient condition for a 
Sobolev map  $f : U \to \r^n$ to satisfy the naturality of the exterior derivative for smooth forms. Our main results are Theorem \ref{th.1} and  Theorem \ref{th.2}. As consequences of these theorem, we can formulate the following special results (corollaries \ref{cor.1} and \ref{cor.2}):

\bigskip

\begin{enumerate}[$\bullet$]
  \item \emph{Let $U$ be a bounded domain in $\r^m$ and  $f \in W^{1,k+1}(U,\r^n)$. Then  the chain rule  (\ref{eq.nat}) holds for any
smooth $k$-forms $\alpha$ on $\r^n$.}
  \item  \emph{Suppose that  $f \in W^{1,k}(U,\r^n)$. If all the $k\times k$ minors of the Jacobian matrix 
  $\left(\frac{\partial f_{\nu}}{\partial x_{\mu}} \right)$ belong to the space 
   $L^{k/(k-1)}(U)$, then  the chain rule (\ref{eq.nat}) holds for any smooth $k$-forms $\alpha$ on $\r^m$. }
\end{enumerate}
 
\bigskip

\textbf{Remarks 1.)}  The first results says in particular that if $f \in W^{1,m}(U,\r^n)$, then the naturality (\ref{eq.nat}) holds for a
smooth form of any degree. See \cite{GTconformal} for more on this case.

\smallskip

\textbf{2.)}  The case $k=n-1$ of the second result has been studied by J. Ball and  V. \v Sver\'ak
\cite{ball,Sv}. In this special case, this result has also been improved in 
by S. M\"uller, T. Qi  and  B.S. Yan. These authors proved in \cite{MQY}  that  this result is also true for $k = n-1$,
$f\in W^{1,n-1}(U;\r^n)$ and  $|\Lambda^k(f)| \in L^q(U)$  for some  $q\geq n/(n-1)$ (instead of $q \geq p/(p-1)$).
See also \cite[page 256]{giaquinta98} for another proof of this result in the context of the theory of Cartesian currents.  

\smallskip

\textbf{3.)}  For convenience, we formulate our results for maps from a bounded domain into euclidean space. However, the 
chain rule  (\ref{eq.nat}) is a local formula and our results also apply to the case of mappings between smooth manifolds.

\section{Measurable differential forms}

Let $U \subset \r^m$ be a domain in $m$-dimensional euclidean space. 
A \emph{measurable differential form} of degree $k$ in $U$ is a measurable function 
$\theta : U \to \Lambda^k (\r^m)$. If $x_1,x_2,\dots,x_m$ is a system of smooth coordinates in $U$, then 
any measurable differential $k$-form writes as
$$
 \theta = \sum_{i_1<i_2<\cdots <i_k}  h_{i_1i_2\dots i_k}(x) 
 dx_{i_1}\wedge dx_{i_2}\wedge dx_{i_k},
$$
where  the coefficients $h_{i_1i_2\dots i_k}$ are measurable functions on $U$.
The form $\theta$ belongs to $L^p(U,\Lambda^k)$ if $h_{i_1i_2\dots i_k}\in L^p(U)$
for all multiindices $i_1i_2\dots i_k$ and similarly 
$\theta\in C^r(U,\Lambda^k)$ if all  $h_{i_1i_2\dots i_k}\in C^r(U)$. If the  coefficients vanish outside a compact subset of $U$, then one writes $\theta\in C_0^r(U,\Lambda^k)$.

\bigskip

Any $k$-form $\theta \in L^p(U,\Lambda^k)$ defines a continuous linear form on the space  $\omega\in C_0^{1}(U,\Lambda^{m-k})$ by the following formula:
$$
 \langle \theta , \omega \rangle = \int_U \theta \wedge \omega.
$$

\bigskip

\begin{definition}\label{def.formweakconv}
A sequence $\{\theta_j \} \subset L^1(U, \Lambda^k)$ is said to 
\emph{converge weakly} to $\theta \in  L^1(U, \Lambda^k)$ if and only if  for every
$\omega \in  C_0^{1}(U, \Lambda^{m-k})$, we have
$$\int_U \theta_j\wedge \omega   \to \int_U \theta\wedge \omega \, .$$
\end{definition}
It is clear that strong convergence in $L^1$ implies weak convergence. The converse is not true.

\begin{definition}
Let $\theta \in L^1_{loc}(U, \Lambda^k)$ be a  $k-$form. If  there exists a $(k+1)-$form
$\psi \in L^1_{loc}(M, \Lambda^{k+1})$ for which the equality
$$\int_U \theta\wedge d\omega = (-1)^{k+1}\int_U\psi\wedge\omega $$
holds for any $\omega \in  C_0^{1}(U, \Lambda^{m-k-1})$, then
 $\psi$ is   called the {\sl weak exterior derivative} of $\theta$
(or the exterior derivative of $\theta$ in the sense of currents) and is denoted 
by $\psi = d\theta$.
The form $\theta \in L^1_{loc}(M, \Lambda^k)$ is \emph{weakly closed} if $d\theta = 0$ in the weak
sense, that is  if 
$$\int_U \theta\wedge d\omega =0 $$
holds for  any  $\omega\in C_0^{1}(U,\Lambda^{m-k-1})$.
\end{definition}

\bigskip

\begin{lemma}
Let $\alpha \in L^1_{loc}(U, \Lambda^k)$ and
$\beta \in  L^1_{loc}(U, \Lambda^{k+1})$.
If there exists a sequence $\{\alpha_j \} \subset  C^1(U, \Lambda^k)$ such
that
$\alpha_j \to \alpha$ and $d\alpha_j \to \beta$ weakly, 
then  $d\alpha = \beta$ in the weak sense.
\end{lemma}

\medskip

\textbf{Proof } For any $\omega \in  C_0^{1}(U, \Lambda^{m-k-1})$, we have
$$\int_U \alpha\wedge d\omega = \lim_{j\to \infty} \int_U \alpha_j\wedge d\omega
= (-1)^{k+1}\int_U d\alpha_j\wedge \omega = 
(-1)^{k+1}\int_U \beta\wedge \omega \, .$$

\qed  


\bigskip

\begin{lemma}\label{lem.Leibniz}
Let $h : U \to \r$ be a bounded function such that $dh \in L^{p'}(\r^n)$ and $\beta \in L^p(U, \Lambda^k)$ such that
$d\beta \in L^{\infty}(U, \Lambda^k)$ where $p' = p/(p-1)$. 
Then $h \cdot \beta  \in L^p(U, \Lambda^k)$ and
\begin{equation}\label{eq,leibniz}
 d(h \cdot \beta) = dh \wedge \beta + h \cdot d\beta. 
\end{equation}
\end{lemma}
 
 \medskip
 
\textbf{Proof}  The equation (\ref{eq,leibniz}) is classic for smooth forms. Use now the density of smooth forms in $L^p$ and
the H\"older inequality to obtain the equality (\ref{eq,leibniz}) in the general case.

\qed

\section{Sobolev mappings}

\textbf{Definition} A map $f : U \to \r^n$ is said to be \emph{bounded} if $f(U) \subset \r^n$ is relatively compact. It belongs to $W^{1,p}(U,\r^n)$ if all its component $(f_1,f_2,\dots,f_n)$ belong to the Sobolev space $W^{1,p}(U,\r)$.

\medskip

Given a map $f \in W^{1,p}(U,\r^n)$, one  defines the pullback of a smooth differential form $\alpha \in C^1(\r^n, \Lambda^k)$ by the following formula:
if 
$$
 \alpha = \sum_{i_1<i_2<\cdots <i_k}  a_{i_1i_2\dots i_k}(y) 
 dy_{i_1}\wedge dy_{i_2}\wedge dy_{i_k},
$$
then
$$
 f^*\alpha = \Lambda^k f (\alpha)  = \sum_{i_1<i_2<\cdots <i_k}  a_{i_1i_2\dots i_k}(f(x)) \; 
 df_{i_1}\wedge df_{i_2}\wedge df_{i_k},
$$
where 
$$
 df_{\nu} = \sum_{\mu = 1}^m \frac{\partial f_{\nu}}{\partial x_{\mu}} dx_{\mu}.
$$
Clearly, $f^*\alpha $ is a differential form with measurable coefficients in $U$ for any 
$\alpha \in C^1(\r^n, \Lambda^k)$.

\bigskip

Let us denote by $Df(x)$ the formal Jacobian matrix of $f$ at the point $x\in U$. This is the $n\times m$ matrix
whose entries are the partial derivatives of $f$:
$$
 Df = \left(\frac{\partial f_{\nu}}{\partial x_{\mu}} \right),
$$
it is defined almost everywhere in $U$.

\bigskip

The pullback operator $ \Lambda^k f$ is represented by the matrix of $k\times k$ minor determinants of $Df(x)$. Indeed we have
 by linear algebra
 $$
 \Lambda^k f (dy_{i_1}\wedge dy_{i_2}\wedge \dots \wedge dy_{i_k}) = \sum_{j_1 < j_2 < \dots < j_k}
 \frac{\partial (f_{i_1}, f_{i_2}, \dots , f_{i_k})}{\partial(x_{j_1},x_{j_2}, \dots,x_{j_k})}
   dx_{j_1}\wedge dx_{j_2}\wedge  \dots \wedge  dx_{j_k} ,
$$
where we have used the old fashioned but convenient  notation 
$
 \frac{\partial (f_{i_1}, f_{i_2}, \dots , f_{i_k})}{\partial(x_{j_1},x_{j_2}, \dots,x_{j_k})} 
$
to denote the entries of the $k\times k$ minor determinant of $Df$. 

\bigskip

We will use the following norm for  $\Lambda^k f$:
$$
 | \Lambda^k f | = \max  \left|\frac{\partial (f_{i_1}, f_{i_2}, \dots , f_{i_k})}{\partial(x_{j_1},x_{j_2}, \dots,x_{j_k})} \right|,
$$
where the max is taken over all ordered $k$-tuple ${i_1}<{i_2}, \dots < {i_k}$;  ${j_1}<{j_2}, \dots < {j_k}$. Observe that
\begin{equation}\label{eq-surlamb}
  |\Lambda f ^k(\alpha)| \leq   | \Lambda^k f | \,  | \alpha | 
\end{equation}
Observe finally that the map $f \mapsto \Lambda^k f$ is non linear for $k\geq 2$.

\section{The class  $\mathcal{F}^k(U,\r^n)$ }

\textbf{Definition} Let us denote by  $\mathcal{F}^k(U,\r^n)$ the class of maps $f : U \to \r^n$
defined as follow:
$$
 f \in \mathcal{F}^k(U,\r^n) \ \Leftrightarrow \ f \in W^{1,1}(U,\r^n) \text{ and } \Lambda^k(f) \in L^1(U).
$$
This definition is motivated by the obvious fact that for  any map  $f \in  \mathcal{F}^k(U,\r^n)$,  the pull back $\alpha \mapsto f^*\alpha$
defines a bounded operator
$$
 \Lambda^k f = f^* : C_0^1(\r^n, \Lambda^k)  \to L^1(U,\Lambda^k).
$$

Observe that the $\mathcal{F}^1(U,\r^n)=W^{1,1}(U,\r^n) $ and that  $\mathcal{F}^k(U,\r^n)$ is not a vector space
for $2 \leq k \leq m$.

\bigskip
 
We denote by $\tau^k$  the initial topology on  $\mathcal{F}^k(U,\r^n)$  induced by the  inclusion $\mathcal{F}^k(U,\r^n) \subset W^{1,1}(U,\r^n)$ and the family of functions 
$$
 \lambda_{\alpha,\omega} : \mathcal{F}^k(U,\r^n) \to \r , \qquad 
  \lambda_{\alpha,\omega}(f) = \int_U f^*\alpha \wedge \omega
$$
where $\alpha \in C^1(\r^n, \Lambda^k)$ and   $\omega \in C^1_0(U, \Lambda^{m-k})$. In other words $\tau^k$ is the  
coarsest topology for which the inclusion $\mathcal{F}^k(U,\r^n) \subset W^{1,1}(U,\r^n)$ is continuous, as well
as all functions $ \lambda_{\alpha,\omega}$.
 
\medskip

Observe that if a sequence $f_j \in \mathcal{F}^k(U,\r^n)$ converges  to a map $f$ in the topology  $\tau^k$, then
$f_k^*\alpha$ converges weakly to $f^*\alpha$ by definition.

\bigskip

An explicit sufficient condition for the  $\tau^k$-convergence in $\mathcal{F}^k(U,\r^n)$ is given in the next result:

\begin{lemma} \label{lem.limit1}
Let ${\{f_j\}} \subset W^{1,1}(U,\r^n)$ be a sequence of mappings which converges 
to a map $f \in  \mathcal{F}^k(U,\r^n)$ in the $W^{1,1}$-topology. Assume that ${\{| \Lambda^k f_j |\}}$ is equi-integrable, i.e.
there exists a function $w \in L^1(U, \r)$ such that $| \Lambda^k f_j |\leq w(x)$ a.e. $x\in U$ for any $j \in \mathbb{N}$.
Then $f_j\to f$in the $\tau^k$ topology.
\end{lemma}

\bigskip

\textbf{Proof} Let  $\alpha \in C^1(\r^n, \Lambda^k)$ be an arbitrary smooth $k$-form on $\r^n$ and $\omega \in C_0^1(U, \Lambda^{m-k})$. 
Since $f_j \to f$ in $W^{1,1}$, we have 
$$
 \lim_{j \to \infty} (f^*_j \alpha) \wedge \omega=  \lim_{j \to \infty} (\Lambda^k f_j)  \alpha\wedge \omega = f^* \alpha\wedge \omega
$$
almost everywhere. Furthermore, we have at every point $x\in U$
$$
 | (f^*_j \alpha)_x\wedge \theta_x | \leq  |\Lambda^k f_j(x)| \,  |\alpha_x | \, |\theta_x| \leq Q \cdot  |\Lambda^k f_j(x)| \leq Q \cdot w(x)
$$
for some constant $Q$. Because $w \in L^1(U, \r)$, the Lebesgue dominated convergence theorem implies that
\begin{equation}\label{limdfd1}
 \lim_{j \to \infty}  \int_U (f^*_j \alpha)\wedge \omega =   \int_U (f^*\alpha)\wedge  \omega.
\end{equation}

\qed

\begin{proposition}
Let  $f\in W^{1,1}(U, \r^n)$ be a map such that 
\begin{enumerate}[a.)]
  \item The $m$-dimensional Hausdorff measure of the image $f(U) \subset \r^n$ is finite;
  \item  $f$ has \emph{essentially finite multiplicity }, i.e. if there exists a constant $Q < \infty$ and a set $E\subset U$
with measure zero such that for every point $y\in \r^n$, 
$$
 \Card \{x \in U \setminus E \tq  f(x) = y \}  \leq Q.
$$
Then $f \in \mathcal{F}^m(U,\r^n)$.
\end{enumerate}
\end{proposition}

\medskip 

This proposition applies e.g. if $f$ is a homeomorphism onto a bounded domain.

\medskip

\textbf{Proof} In that case, $|\Lambda^k(f)|$ belongs to $L^1$  by the area formula
(see e.g. \cite[page 220]{giaquinta98}).

\qed

\bigskip

\textbf{Remark}  In \cite[page 229]{giaquinta98}, Giaquinta introduce a class of maps $\mathcal{A}_1(U,\r^n)$ which is
very similar to our class $\mathcal{F}^m(U,\r^n)$ (where $m= \dim (U)$). The main difference is that the condition 
$f\in W^{1,1}(U, \r^n)$ is relaxed to the assumption that $f$ is approximately differentiable almost everywhere. In any
case, we have a continuous embedding
$$
 W^{1,1}(U, \r^n) \subset \mathcal{A}_1(U,\r^n).
$$

\section{$k$-stable maps in $\mathcal{F}^k(U,\r^n)$ }

\textbf{Definition} A map $f \in \mathcal{F}^k(U,\r^n)$ is said to be $k$-\emph{stable} 
if it belongs to the  closure of $C^1(U,\r^n)$ in the $\tau^k$ topology, i.e. there exists a sequence 
of smooth maps  converging to $f$ in the $\tau^k$ topology.
We denote by $\mathcal{S}^k(U,\r^n) \subset \mathcal{F}^k(U,\r^n)$ the  set of $k$-stable maps :
$$
 \mathcal{S}^k(U,\r^n) = \overline{C^1(U,\r^n)}^{\, \tau^k}  \subset \mathcal{F}^k(U,\r^n).
$$
observe that $W^{1,k}(U, \r^n) \subset \mathcal{S}^k(U,\r^n)$.

\medskip

The pullback of a closed form by a stable map is again a closed form:

\begin{proposition} \label{pullbackclosed}
Let $f \in \mathcal{S}^k(U,\r^n)$  be $k$-stable map and  $\alpha \in C^1(\r^n, \Lambda^k)$.
If $\alpha$ is closed, then $f^*\alpha$ is weakly closed.
\end{proposition}

\medskip

\textbf{Proof} Because $f \in \mathcal{S}^k(U,\r^n)$, there exists a  sequence $\{f_j\}$ of smooth maps
converging to $f$ in the $\tau^k$-topology.  
Assume   that $d\alpha = 0$, then for any 
$\phi \in C_0^1(U, \Lambda^{m-k-1})$ we have 
$$
 \int_U (f^*_j \alpha)\wedge d\phi  =   (-1)^{k+1} \int_U  d(f^*_j \alpha)\wedge \phi  
=  (-1)^{k+1}   \int_U f^*_j (d\alpha)\wedge \phi  = 0.
$$
We  thus have 
\begin{equation}\label{}
 \int_U (f^*\alpha)\wedge d\phi =
 \lim_{j \to \infty}  \int_U (f^*_j \alpha)\wedge d\phi  =  0,
\end{equation}
for any $\phi \in C_0^1(U, \Lambda^{m-k-1})$. This means that  $f^*\alpha$ is weakly closed.

\qed

\bigskip

\begin{proposition}\label{prop.critkstable}
Let $f \in W^{1,1}(U,\r^n)$  be a map such that 
$$
 \inf_{\{f_j\}} \int_U \left( \sup_j |\Lambda^k f_j|\right) \, dx < \infty,
$$
where the infimum is taken over the set of all sequences ${\{f_j\}}$ of smooth maps 
such that \ $ \|f_j - f \|_{W^{1,1}} \to 0$. Then $f \in \mathcal{S}^k(U,\r^n)$.
\end{proposition}

\medskip

\textbf{Proof}  By mollification, we know that the set   sequences ${\{f_j\}}$ of smooth maps 
such that \ $ \|f_j - f \|_{W^{1,1}} \to 0$ is not empty. We can then apply Lemma \ref{lem.limit1}.

\qed

\bigskip

\section{$k^\dagger$-stable maps}

\begin{definition}
We define the space $\mathcal{F}^{k^\dagger}(U,\r^n)$ by
$$
  \mathcal{F}^{k^\dagger}(U,\r^n)  = 
  \begin{cases}
  \mathcal{F}^n(U,\r^n)  & \text{if $k=n$} , \\  \\
  \mathcal{F}^k(U,\r^n) \cap  \mathcal{F}^{k+1}(U,\r^n) & \text{if $0 \leq k<n$}. \\
\end{cases}
$$
The $\tau^{k^\dagger}$ topology is defined  for $k<n$ to be the initial topoply for which both inclusions
$$
  \mathcal{F}^{k^\dagger}(U,\r^n)  \subset  \mathcal{F}^{k}(U,\r^n) 
  \qquad \text{and} \qquad
   \mathcal{F}^{k^\dagger}(U,\r^n)  \subset  \mathcal{F}^{k+1}(U,\r^n) 
$$
are continuous. For $k=n$, we simply define $\tau^{k^\dagger} = \tau^k$.
\end{definition}

\medskip

We then say that a map $f : U \to \r^n$ is $k^\dagger$\emph{-stable}  if it belongs to the  closure of $C^1(U,\r^n)$ in 
the space  $\mathcal{F}^{k^\dagger}(U,\r^n)$ for the $\tau^{k^\dagger}$ topology.

\bigskip

Observe the following elementary

\begin{lemma}
 A  map $f : U \to \r^n$ is $k^\dagger$-stable if and only if there exists a sequence $\{ f_j \} \subset C^1(U,\r^n)$   
 of smooth maps which  weakly converges to $f$ in both spaces $\mathcal{F}^{k}(U,\r^n)$ and $\mathcal{F}^{k+1}(U,\r^n)$. 
\end{lemma}

\medskip
\qed

\bigskip

\begin{proposition}\label{prop.critk+stable}
Let $f \in W^{1,1}(U,\r^n)$ be a map such that for some $k<n$,
$$
 \inf_{\{f_j\}} \int_U \left( \sup_j (|\Lambda^k f_j| +  |\Lambda^{k+1} f_k|) \right) \, dx < \infty,
$$
where the infimum is taken over all sequences ${\{f_j\}}$ of smooth maps such that \ $ \|f_j - f \|_{W^{1,1}} \to 0$. Then $f \in \mathcal{S}^{k^\dagger}(U,\r^n)$.
\end{proposition}

\textbf{Proof} This follows directly from Proposition \ref{prop.critkstable} and the previous lemma.

\qed

One can rephrase this Proposition as follow.  Let $f\in W^{1,1}(U,\r^n)$, and assume that there exists 
 a sequence of    smooth maps$\{ f_j \} \subset C^1(U,\r^n)$ and  a function $w \in L^1(U, \r)$ such that 
$$|\Lambda^k f_j(x)| + |\Lambda^{k+1} f_j(x)|\leq w(x)$$
a.e. $x\in U$ for any $j \in \mathbb{N}$. Then  $f$ is  $k^\dagger$-stable.

\bigskip

The naturality of the exterior differential holds for  $k^\dagger$-stable maps:

\begin{theorem}\label{th.1}
Let $f \in  \mathcal{S}^{k^\dagger}(U,\r^n)$  be $k^\dagger$-stable map, and let 
$\alpha \in C^1(\r^n, \Lambda^k)$ be a smooth $k$-form in $\r^m$, then $f^*\alpha \in L^1(U, \Lambda^k)$, $f^*d\alpha \in L^1(U, \Lambda^{k+1})$ and the equation
$$
 df^*\alpha  = f^*d\alpha 
$$
holds in the weak sense.
\end{theorem}

\bigskip

\textbf{Proof} By hypothesis, there exists a sequence of smooth mappings $f_j \in C^1(U,\r^n)$
which converges to $f$ in   $\mathcal{F}^k(U,\r^n)$ and $\mathcal{F}^{k+1}(U,\r^n)$
for both the $\tau^k$ and $\tau^{k+1}$ topologies.

Let  $\alpha \in C^1(\r^n, \Lambda^k)$ be an arbitrary smooth $k$-form on $\r^m$ and $\theta  \in C_0^1(U, \Lambda^{m-k})$.
 By hypothesis, we have 
\begin{equation}\label{limdfd1a}
 \lim_{j \to \infty}  \int_U (f^*_j \alpha)\wedge \theta =   \int_U (f^*\alpha)\wedge \theta.
\end{equation}
We also have
\begin{equation}\label{limdfd2}
 \lim_{j \to \infty}  \int_U (f^*_j \beta)\wedge \phi=   \int_U (f^*\beta)\wedge \phi
\end{equation}
for any  $\beta\in C^1(\r^n, \Lambda^{k+1})$ and $\phi \in C_0^1(U, \Lambda^{m-k-1})$.

Let us now choose $\beta = d\alpha$ and $\theta = d\phi$, we then have $df_j^*\alpha = f_j^*d\alpha$
for any $j\in \mathbb{N}$ because both $\alpha$ and $f_j$ are of class $C^1$, this imples that
$$
   \int_U (f_j^*\alpha)\wedge d\phi = (-1)^{k+1}\int_U  d(f_j^*\alpha)\wedge \phi 
   = (-1)^{k+1}\int_U (f_j^*d\alpha)\wedge \phi.
$$
Applying (\ref{limdfd1}) and (\ref{limdfd2}) one gets  then
\begin{align*}
   \int_U (f^*\alpha)\wedge d\phi & =  \lim_{j \to \infty}  \int_U (f^*_j \alpha)\wedge d\phi 
   \\ &   =   \lim_{j \to \infty}  (-1)^k\int_U  f_j^*(d\alpha)\wedge \phi 
   \\ &   =   (-1)^k\int_U  f^*(d\alpha)\wedge \phi 
\end{align*}
for any $\phi \in C_0^1(U, \Lambda^{n-k-1})$, this means precisely that $d(f^*\alpha) = f^*(d\alpha)$ in the weak sense.

\qed

\bigskip

\begin{corollary}\label{cor.1}
Let $U$ be a domain in $\r^m$ and  $f \in W^{1,k+1}(U,\r^n)$. Then  the naturality (\ref{eq.nat}) holds for any
smooth $k$-forms $\alpha$ on $\r^n$.
\end{corollary}

\textbf{Proof} This follows from the fact that  $W^{1,k+1}(U,\r^n) \subset  \mathcal{S}^{k^\dagger}(U,\r^n)$.

\qed

\section{Another class of maps}

We denote by $\mathcal{S}^{k}_{q,p}(U,\r^n)$ the class of maps $f \in \mathcal{S}^{k}(U,\r^n)$
such that 
$$
 |df| \in L^p(U)  \qquad \text{and}  \qquad  |\Lambda^k(f)| \in L^q(U).
$$
Observe that  $\mathcal{S}_{q,p}^k(U,\r^n) \subset W^{1,p}(U,\r^n)$.

\bigskip

\begin{theorem}\label{th.2}
Let $f \in \mathcal{S}_{q,p}^k(U,\r^n)$,  and assume $1\leq p \leq \infty$, $q = p/(p-1)$.

Let $\alpha \in C^1(\r^n, \Lambda^k)$ be a smooth $k$-form in $\r^n$, then $f^*\alpha \in L^1(U, \Lambda^k)$, 
$f^*d\alpha \in L^1(U, \Lambda^{k+1})$ and the chain rule
$$
 df^*\alpha  = f^*d\alpha 
$$
holds in the weak sense.
\end{theorem}

\bigskip

\textbf{Proof} 
Observe that $f^*\gamma$ is weakly closed for any
closed  $k$-form $\gamma  \in C^1(\r^m, \Lambda^k)$  by Proposition \ref{pullbackclosed}.

Suppose  first that $\alpha = a \cdot \gamma$ where  $\gamma  \in C^1(\r^m, \Lambda^k)$ is a closed $k-$form  and 
that $a\in C^1(\r^n)$ is a function. Then $f^*a = a \circ f \in W^{1,1}(U)$ and $df^*a = f^* da$ (see e.g. \cite[Theorem 7.8]{gilbarg}).
Because $f \in \mathcal{S}_{q,p}^k(U,\r^n)$, we have in fact $|df^*a| \in L^p(U)$ and $|f^*(\gamma)| \leq \ |\Lambda^{k+1} f_j(x)|  \cdot |\gamma| \in L^{p'}(U)$.
Therefore we have by Lemma \ref{lem.Leibniz}:
\begin{align*}
d f^*\alpha & = d f^*( a \cdot \gamma )
       \\ & =    d ( f^*a \cdot f^*\gamma )
        \\ & =  d ( f^*a) \wedge f^*\gamma  +   ( f^*a) \cdot \underbrace{(df^*\gamma )}_{=0}
        \\ & =  d ( f^*a) \wedge f^*\gamma 
        \\ & =   ( f^*da) \wedge f^*\gamma 
        \\ & =   f^*(da \wedge \gamma)  
         \\ & =    f^*(d \alpha)      
\end{align*}

Consider now an arbitrary smooth $k$-form on $\r^n$. It can be written as a sum
$$
 \alpha = \sum_{i_1<i_2<...<i_k} a_{i_1i_2...i_k}(x) \; dy_{i_1}\wedge dy_{i_2}\wedge \cdots dy_{i_k},
$$
where $ a_{i_1i_2...i_k}(x) $ is an element in $C^1(\r^n)$. Since $dy_{i_1}\wedge dy_{i_2}\wedge \cdots dy_{i_k}$
is a closed (in fact exact) form, the proof is complete.

\qed

\bigskip

\begin{corollary}\label{cor.2}  
 Suppose that  $f \in W^{1,k}(U,\r^m)$ and   $\Lambda^k(f) \in L^{k/(k-1)}(U)$, then  the naturality (\ref{eq.nat}) holds for any
smooth $k$-forms $\alpha$ on $\r^m$.\end{corollary}

\medskip

\textbf{Proof.}  The hypothesis imply that  $f \in \mathcal{S}_{q,p}^k(U,\r^n)$ .

\qed

\vfill

\makeatother

{Vladimir Gol'd'shtein} \\
{Department of Mathematics,  
Ben Gurion University of the Negev,  \\
P.O.Box 653, Beer Sheva, Israel}  \\
{vladimir@bgumail.bgu.ac.il}

\medskip

{Marc Troyanov} \\ 
Institut de G{\'e}om{\'e}trie, alg{\`e}bre et topologie (IGAT)
B{\^a}timent BCH \\
 \'Ecole Polytechnique F{\'e}derale de
Lausanne, 1015 Lausanne - Switzerland \\
{marc.troyanov@epfl.ch}

\end{document}